\documentclass[10pt]{amsart}
\usepackage{amssymb, enumitem}
\usepackage[all]{xy}
\usepackage{tikz}
\usepackage{relsize}
\usepackage{aliascnt}
\usepackage{verbatim}
\usepackage[english]{babel}
\usepackage{enumitem}
\usepackage{mathtools}

	\usepackage{amsmath} %formulars and symbols
	\usepackage{amssymb} %more symbols
	\usepackage{mathrsfs}  %script letters
	\usepackage{mathtools} %curved arrows
	\usepackage[T1]{fontenc} %if a font-nerd reads this
	\usepackage[allcolors=black]{hyperref} %clickable links
	\usepackage[utf8]{inputenc} %some font stuff??
	\usepackage{bm} %bold symbols

	\usepackage{graphicx} %graphics
	\usepackage{tikz-cd} %commutative diagrams
	\usepackage[noadjust]{cite} %bibliography, noadjust= no silly spaces before citation
	\usepackage{todonotes} %insert todo notes in text
	\usepackage{setspace} %Zeilenabstand Titelseite

	\usepackage{amsthm} %loads the theorem package

\def\today{\number\day\space\ifcase\month\or   January\or February\or
	March\or April\or May\or June\or   July\or August\or September\or
	October\or November\or December\fi\   \number\year}

\theoremstyle{plain}

\newaliascnt{thmCt}{lma}

\aliascntresetthe{thmCt}

\newaliascnt{corCt}{lma}

\aliascntresetthe{corCt}

\newaliascnt{propCt}{lma}

\aliascntresetthe{propCt}

\newtheorem*{thm*}{Theorem}
\newtheorem*{cor*}{Corollary}
\newtheorem*{prop*}{Proposition}

\theoremstyle{plain}
\newtheorem{thmintro}{Theorem} %Use this counter for things in the introduction

\newaliascnt{introcorCt}{thmintro}

\aliascntresetthe{introcorCt}

\newaliascnt{egintroCt}{thmintro}
\newtheorem{egintro}[egintroCt]{Example}
\newaliascnt{propintroCt}{thmintro}
\newtheorem{propintro}[propintroCt]{Proposition}
\theoremstyle{definition}

\newaliascnt{introquesCt}{thmintro}

\aliascntresetthe{introquesCt}

\newaliascnt{pbmintroCt}{thmintro}

\aliascntresetthe{introquesCt}

\newaliascnt{dfnintroCt}{thmintro}

\aliascntresetthe{introquesCt}

\newaliascnt{pgrCt}{lma}

\aliascntresetthe{pgrCt}

\newaliascnt{dfCt}{lma}
\newtheorem{df}[dfCt]{Definition}
\aliascntresetthe{dfCt}

\newtheorem*{df*}{Definition}

\newaliascnt{remCt}{lma}

\aliascntresetthe{remCt}

\newaliascnt{remsCt}{lma}

\aliascntresetthe{remsCt}

\newaliascnt{egCt}{lma}

\aliascntresetthe{egCt}

\newaliascnt{egsCt}{lma}

\aliascntresetthe{egsCt}

\newaliascnt{qstCt}{lma}

\aliascntresetthe{qstCt}

\newaliascnt{pbmCt}{lma}

\aliascntresetthe{pbmCt}

\newaliascnt{notaCt}{lma}

\aliascntresetthe{notaCt}

\theoremstyle{thmCt}

\newtheorem*{conj*}{Conjecture}

\makeatletter
\@namedef{subjclassname@2020}{\textup{2020} Mathematics Subject Classification}
\makeatother
	
		\newcommand{\N}{\mathbb N}

		\newcommand{\C}{\mathbb C}

		\DeclareMathOperator{\Aut}{Aut}

		\DeclareMathOperator{\id}{id}

		\DeclareMathOperator{\Ad}{Ad}
		
		\DeclareMathOperator{\ev}{ev}

		\DeclareMathOperator{\Prob}{Prob}

\title{Actions on classifiable C*-algebras without equivariant property (SI)}

\author[E. Gardella]{Eusebio Gardella} 
\address[Eusebio Gardella]{Department of Mathematical Sciences,
	Chalmers University of Technology and University of Gothenburg,
	Gothenburg,
	SE-412 96,
	Sweden.} 
\email{gardella@chalmers.se}
\urladdr{http://www.math.chalmers.se/~gardella/}

	\author[J. Kranz]{Julian Kranz}
\address[Julian, Kranz]{Data Science: Machine Learning and Data Engineering,
	Prof. Dr. Fabian Gieseke,
	Leonardo-Campus 3,
	48149 M\"unster,
	Deutschland}
\email{julian.kranz@uni-muenster.de}
\urladdr{https://sites.google.com/view/juliankranz/}

\author[A. Vaccaro]{Andrea Vaccaro}
\address[Andrea Vaccaro]{Mathematisches Institut, Fachbereich Mathematik und Informatik der
Universit\"at M\"unster, Einsteinstrasse 62, 48149 M\"unster, Germany.
}
\email{avaccaro@uni-muenster.de}
\urladdr{https://sites.google.com/view/avaccaro/home}

\thanks{
EG was supported by a research grant of 
the Swedish Research Council.
JK was partially supported by the Engineering and Physical Sciences Research Council [Grant Ref: EP/X026647/1] and by the Deutsche Forschungsgemeinschaft (DFG, German Research Foundation) - Project-ID 427320536 - SFB 1442, as well as by Germany's Excellence Strategy EXC
	2044 390685587, Mathematics Münster: Dynamics-Geometry-Structure.
	AV was supported by the Deutsche Forschungsgemeinschaft (DFG, German Research Foundation) under Germany’s Excellence Strategy EXC 2044 –390685587, Mathematics M\"unster: Dynamics–Geometry–Structure, through SFB 1442, by the ERC Advanced Grant 834267 - AMAREC}


\begin{document}

\begin{abstract}
We exhibit examples of actions of countable discrete groups on both simple and non-simple nuclear stably finite C*-algebras that are tracially amenable but not amenable. We furthermore obtain that, under the additional assumption of strict comparison, amenability is equivalent to tracial amenability plus the equivariant analogue of Matui--Sato's property (SI). By virtue of this equivalence,
our construction yields the first known examples of 
actions on classifiable C*-algebras 
that do not have equivariant property (SI). 
We moreover show that such actions can be chosen to absorb the trivial action on the universal UHF algebra, thus proving
that equivariant $\mathcal{Z}$-stability does not in general 
imply equivariant property (SI).
\end{abstract}
\maketitle

%\setcounter{tocdepth}{2}
%\tableofcontents

\section{Introduction}
Property (SI) was isolated in \cite{MatuiSatoActa} by Matui and Sato in the context of the Toms--Winter regularity conjecture. For nuclear C*-algebras it is a consequence of strict comparison, and since its introduction
it has proved to be a powerful tool for converting traits and features of tracial ultrapowers, or more abstractly properties involving tracial 2-norms, into attributes of C*-norm ultrapowers and properties expressible with the C*-norm. One of the most notorious example is Matui--Sato's work itself \cite{MatuiSatoActa}, and its later generalizations \cite{KirRor_Central_2014,TWW, sato:boundary},
proving that, under some assumptions on the trace space, simple separable nuclear unital infinite-dimensional C*-algebras with strict comparison are $\mathcal{Z}$-stable.

The techniques developed in the context of the Toms--Winter conjecture that we alluded to above were also profoundly influential in the study of actions of amenable groups on \emph{classifiable} C*-algebras. This term denotes the class of all separable nuclear simple $\mathcal{Z}$-stable C*-algebras that satisfy the Universal Coefficient Theorem and which, by now, are known to be classifiable by means of $K$-theoretic and tracial invariants (see \cite{WinterICM,WhiteICM} for an overview on the topic). For actions on such algebras, \emph{equivariant $\mathcal{Z}$-stability}, namely the property of absorbing up to cocycle conjugacy
the identity action on the Jiang--Su algebra $\mathcal{Z}$, is a powerful and desirable regularity feature that serves a purpose similar to that of $\mathcal{Z}$-stability in classification. Among its numerous applications, this notion  has recently seen a great deal of interest due to the role it plays in establishing $\mathcal{Z}$-stability, and ultimately classifiability, of crossed products originating from actions on classifiable C*-algebras (see \cite{sato,GaborSI,GHV,Wouters,ESPA,ESPJA}).

Most of the proofs that verify equivariant $\mathcal{Z}$-stability for actions on stably finite C*-algebras follow the same pattern, inspired by the work in \cite{MatuiSatoActa}: first an equivariant McDuff-type condition is established in the tracial ultrapower of the algebra acted upon. Then this property is transferred to the C*-norm ultrapower, where equivariant $\mathcal{Z}$-stability is verified using an equivalent formulation stating the existence of an invariant copy of $\mathcal{Z}$ in the central sequence algebra. The tool allowing the alluded transfer is an equivariant version of Matui--Sato's property (SI), 
implicitly appearing already in \cite{MS:1, MS:2, sato}, which is now referred to as \emph{equivariant property (SI)} (\cite[Definition 2.7]{GaborSI}; see Definition \ref{df-SI}).

Equivariant property (SI) is known to automatically hold for actions of amenable groups on infinite-dimensional simple separable nuclear C*-algebras with strict comparison (\cite{GaborSI}), a result which has recently been extended also to amenable actions (in the sense of \cite{anantharaman1987amenableC}) of non-amenable groups (\cite{GaborSIamenable}). We show here that non-amenable actions may fail to have equivariant property (SI), even when equivariantly $\mathcal{Z}$-stable.

\begin{egintro}\label{eg:SI}
	Let $F_n$ be the free group on $n\geq 2$ generators. There exists an action $\alpha\colon F_n\to \Aut(A)$ on a C*-algebra $A$ with the following properties:
	\begin{enumerate}
		\item $A$ is a unital simple AF-algebra,
		\item \label{item2} $\alpha$ is tracially amenable but does not have the weak containment property, in particular it is not amenable,
		\item $\alpha$ is conjugate to its tensor product with the identity on the universal UHF algebra. In particular, $\alpha$ is equivariantly $\mathcal Z$-stable. 
		\item \label{item4} $\alpha$ does not have equivariant property (SI).
	\end{enumerate}
\end{egintro}

The action constructed in \autoref{eg:SI} is the first 
known example of an equivariantly $\mathcal{Z}$-stable
action of a discrete group on a classifiable C*-algebra without equivariant 
property (SI).
We remark that, in the non-equivariant setting, property (SI) is an automatic consequence of $\mathcal{Z}$-stability, if nuclearity is assumed (\cite{rordamStable, MatuiSatoActa}).
Our \autoref{eg:SI} shows the failure of the analogous implication in the equivariant setting, answering negatively a question of Mikael R{\o}rdam. In hindsight, this should not come
as a complete surprise: as we mentioned, the arguments showing that $\mathcal{Z}$-stability implies property (SI) crucially rely on nuclearity. Amenability often acts as a natural counterpart of nuclearity in the equivariant setting, so it makes sense that the lack of it constitutes an obstacle when trying to generalize this implication to actions.

As hinted by item \eqref{item2}, \autoref{eg:SI} originates from the comparison between the notions of \emph{amenability} and \emph{tracial amenability} for actions on C*-algebras. The former was introduced in \cite{anantharaman1987amenableC} as a C*-algebraic adaptation of the definition of amenability on von Neumann algebras. 
\begin{df}[Amenability]
An action $\alpha \colon G \to \text{Aut}(A)$ of a discrete group on a C*-algebra
is \emph{amenable} if the induced action $\alpha^{**}$ on the double dual $A^{**}$ is a \emph{von Neumann amenable} in the sense of \cite{amenableVNalgebra}.
\end{df}

There has recently been a renewed interest in this property, leading to some significant progress in the study of C*-dynamics, mainly in the form of generalizations of techniques that, up until few years ago, were limited to actions of amenable groups; see for instance \cite{BusEchWil22, OzawaSuzuki} and the references therein. 
On the other hand, tracial amenability was defined in \cite{ESPJA} as a tracial analogue of amenability, which is weak enough to be verified in concrete examples, yet strong enough to guarantee interesting structural properties for the crossed product (such as pure infiniteness). We recall its definition below.

\begin{df}[Tracial Amenability]
An action $\alpha \colon G \to \text{Aut}(A)$ of a discrete group on a stably finite unital C*-algebra is said to be \emph{tracially amenable} if the induced action $\alpha^{**}_{\text{fin}}$ on the \emph{finite part} of the double dual $A^{**}_{\text{fin}}$ is von Neumann amenable.
\end{df}

The main motivation for introducing this definition is that it is currently not known whether amenable actions of non-amenable groups on simple unital stably finite C*-algebras exist. Purely infinite examples have been found in \cite{SuzukiEquivariant, OzawaSuzuki}, as well as simple stably finite \emph{non-unital} ones in \cite{suzuki1, suzuki2}. More recently, in  \cite{ suzuki3} Suzuki also built amenable actions on simple unital C*-algebras which are neither stably finite nor purely infinite. On the other hand, in \cite{ESPJA} it is shown that an action is tracially amenable if and only if the induced action on the trace space is topologically amenable, allowing to use techniques developed within the Elliott Classification Program to construct concrete examples.

While it is clear that amenability implies tracial amenability, in \cite{ESPJA} it was left open whether the converse holds. The actions in \autoref{eg:SI} were originally built to show that these two conditions are different. Indeed, our construction is based on the tracially amenable actions built in \cite[Example 2.11]{ESPJA} and uses the main result of \cite{KOS} to produce an approximately inner perturbation of such actions with an invariant state\footnote{A similar remark can also be found after \cite[Lemma 2.3]{suzukitight}.}, which thus cannot be amenable by \cite[Proposition 3.5]{OzawaSuzuki}. As pointed out to us by Damian Ferraro, having an invariant state implies furthermore the failure of the {weak containment property}, since this condition forces the stabilizers of states to be amenable. We recall that an action $\alpha \colon G \to \Aut(A)$ has the \emph{weak containment property} if the canonical surjection $A\rtimes_{\mathrm{max}} G \to A \rtimes_{\mathrm{red}} G$ is an isomorphism, and that this property is automatic for amenable actions of discrete groups by  \cite[Proposition 4.8]{anantharaman1987amenableC}.

Item \eqref{item4} of \autoref{eg:SI} is a consequence of the following equivalence, showing that the difference between amenability and tracial amenability is precisely equivariant property (SI), which once again fills the gap between the tracial and the operator-norm world. We prove moreover that the discrepancy between the two properties can also be identified in the notion of \emph{separable stability} of the fixed point-algebra of the intersection of the central sequence algebra with the trace kernel ideal. We defer the precise definitions to Section \ref{S.2}.
As we mentioned previously, the implication \eqref{SI1}$\Rightarrow$\eqref{SI2} was obtained in \cite{GaborSIamenable} (a partial version of this equivalence has also been independently proved in JK's PhD thesis \cite{JulianPhD}).

\begin{thmintro}\label{thm:TrAmSI}
	Let $\alpha\colon G\to\Aut(A)$ be an action of a countable discrete group $G$ on a unital simple separable nuclear stably finite C*-algebra $A$ with strict comparison. The following are equivalent:
	\begin{enumerate}
		\item\label{SI1} $\alpha$ is amenable,
		\item\label{SI2} $\alpha$ is tracially amenable and has equivariant property (SI).
		\item\label{SI3} $\alpha$ is tracially amenable and $(J_{T(A)} \cap A')^{\alpha_\omega}$ is separably stable.
	\end{enumerate}
\end{thmintro}

We furthermore build tracially amenable actions which are not amenable in the non-simple setting, more precisely on the Toeplitz algebra $\mathcal{T}$.

\begin{egintro}\label{eg-Toeplitz}
	For any $n\geq 2$, there is a tracially amenable, non-amenable action $\alpha\colon F_n\to \Aut(\mathcal{T})$ of the free group on $n$ generators on the Toeplitz algebra.
\end{egintro}

As we will see, the reason why amenability fails in this case is that the action leaves the ideal of compact operators invariant.

%This illustrates that the definition of tracial amenability given in \cite{ESPJA} might not be the correct one for actions on non-simple C*-algebras and that an appropriate definition should take into account tracial weights on all ideals. 

It remains an open problem (see \cite[Problem~D]{ESPJA})
to determine whether there exist amenable actions of nonamenable groups on simple unital
stably finite C*-algebras.
The proof of \autoref{eg:SI} shows that such actions cannot be obtained \emph{abstractly} and directly 
from the classification theorem,
since in this context amenability is not determined by the induced action on the Elliott invariant. To further illustrate the difficulties surrounding this problem, 
we show another obstacle in obtaining such actions by \emph{explicit} 
constructions using amenability of the action of $F_n$ 
on its Gromov boundary $X=\partial F_n$ or any other 
boundary action in the sense of 
\cite[Definition 3.8]{Kal_Ken_2017}.

\begin{propintro}\label{thm-noembedding}
	Let $\alpha\colon G\to \Aut(A)$ be an action of a non-amenable discrete group on a unital, simple C*-algebra with $T(A)\neq \emptyset$. Let $X\not=\{\ast\}$ be a $G$-boundary. Then there is no $G$-equivariant unital completely positive map $C(X)\to A$. 
\end{propintro}

In particular, amenable actions of non-amenable groups on simple, stably finite, unital C*-algebras cannot be constructed as equivariant inductive limits of the form $\varinjlim C(X)\otimes A_n$ for a $G$-boundary $X$. Since most explicit examples of amenable actions of non-amenable groups on compact spaces arise as boundary actions, this may help explain why the problem of finding such actions is so difficult to tackle at the moment.

\subsection*{Acknowledgements}
We thank Jamie Gabe and Stuart White for instructive feedback. We thank furthermore the anonymous referee for indicating to us the equivalence with item \eqref{SI3} in \autoref{thm:TrAmSI}. JK thanks
Mikael R{\o}rdam for valuable conversations, and
AV thanks G\'abor Szab\'o for helpful discussions. We are moreover grateful to Alcides Buss and Damian Ferraro for pointing out to us that Example~\ref{eg:SI} fails to satisfy the weak containment property.

%JK thanks Mikael R{\o}rdam 
%for asking whether there could exist
%an equivariantly
%$\mathcal{Z}$-stable action without property (SI), and AV thanks G\'abor Szab\'o for interesting and helpful discussions.
Some 
of the results in this paper have already appeared in 
JK's PhD thesis \cite{JulianPhD}.

\section{Preliminaries} \label{S.2}

%Let $\alpha\colon G\to \Aut(A)$ be an action of a discrete group $G$ on a C*-algebra $A$. Then $\alpha$ extends uniquely to an action $\alpha^{**}\colon G\to \Aut(A^{**})$ on the enveloping von Neumann algebra of $A$. We equip $\ell^\infty(G)\overline \otimes A^{**}$ with the diagonal action$\lambda\otimes \alpha^{**}$ where $\lambda\colon G\to \Aut(\ell^\infty(G))$ denote the left translation action. We moreover identify $A^{**}$ with the subalgebra $1\otimes A^{**}\subseteq \ell^\infty(G)\overline \otimes A^{**}$. 
%
%\begin{df}[\cite{anantharaman1987amenableC}]
%	An action $\alpha\colon G\to \Aut(A^{**})$ of a discrete group $G$ on a C*-algebra $A$ is called \emph{amenable} if there is a $G$-equivariant conditional expectation 
%	\[\ell^\infty(G)\overline \otimes A^{**}\to A^{**}.\]
%\end{df}
%An action of a discrete group $G$ on a compact space $X$ is called amenable, if the induced action of $G$ on $C(X)$ is amenable in the above sense. 

Fix a free ultrafilter $\omega\in \beta\N\setminus \N$ and let $A$ be a separable unital C*-algebra. The \emph{(C*-norm) ultrapower} is the C*-algebra
\[
A_\omega = \ell^\infty(A) / \{ (a_n)_{n \in \N} \in \ell^\infty(A) \colon \lim_{n \to \omega} \| a_n \| = 0 \}.
\]
We identify $A$ with the subalgebra of all (classes of) constant sequences in $A_\omega$ and we denote by $A_\omega\cap A'$ its relative commutant.

Let $T(A)$ be the set of all tracial states of $A$, which we shall simply refer to as \emph{traces}. Given a non-empty
subset $X \subseteq T(A)$, we define the 2-seminorm
\[
\| a \|_{2, X} = \sup_{\tau \in X} \tau(a^*a)^{1/2},
\]
for all $a\in A$.
A \emph{limit trace} on $A_\omega$ is a trace $\tau \in T(A_\omega)$ for which there is a sequence $(\tau_n)_{n \in \N} $ in $T(A)$ satisfying
\[
\tau((a_n)_{n \in \omega}) = \lim_{n \to \omega} \tau_n(a_n),
\]
for all $(a_n)_{n \in \N} \in A_\omega$.
We denote by $T_\omega(A)$ the set of limit traces on $A_\omega$. The \emph{trace kernel ideal} $J_{T(A)}$ is defined
as 
\[J_{T(A)}=\{a \in A_\omega \colon \| a \|_{2, T_\omega(A)} = 0\}.\]
The \emph{tracial ultrapower} of $A$ is defined to be the quotient $A^\omega = A_\omega/J_{T(A)}$. When $\| \cdot \|_{2, T(A)}$ is a norm on $A$ (for instance when $A$ is simple), then the map sending each element in $A$ to the corresponding constant sequence into $A^\omega$ is an embedding. In such case, we denote by $A^\omega \cap A'$ the relative commutant of the image of this embedding. It is well-known that the restriction of the quotient map $q \colon A_\omega \to A^\omega$ to $A_\omega \cap A'$ is surjective onto $A^\omega \cap A'$ (see for instance \cite[Theorem 3.3]{KirRor_Central_2014}).

Let $G$ be a countable discrete group, and fix an action $\alpha \colon G \to \Aut(A)$. Acting coordinate-wise, $\alpha$ naturally induces actions $\alpha_\omega$ and $\alpha^\omega$ on $A_\omega$ and $A^\omega$ respectively. Given $B \subseteq A$, we let $B^{\alpha}=\{b\in B\colon \alpha_g(b)=b \mbox{ for all } g\in G\}$ be the corresponding fixed point algebra.

Let $\mathcal{K}$ denote the algebra of compact operators on $\ell^2$. 
A C*-algebra $A$ is said to be \emph{separably stable} if for every separable subalgebra $B \subseteq A$ there exists a separable C*-algebra $C \subseteq A$ containing $B$ such that $C \otimes \mathcal{K}(\ell^2) \cong C$. In \cite[Proposition 6.10]{CGSTW} separable stability of $A$ is proved to be equivalent to the \emph{Hjelmborg--R{\o}rdam criterion}, stating that if a positive contraction $e \in A_+$ and $\varepsilon> 0$ are given, then there is a contraction $s \in A$ such that $\| e - s^*s \| < \varepsilon$ and $\| (s^*s)(ss^*) \| < \varepsilon$. This equivalent formulation will be used in the proof of \autoref{thm:TrAmSI}.

\section{Results and Examples}
In this section, we will exploit different equivalent reformulations of amenability and tracial amenability to obtain our main results and examples. We refer to \cite[Theorem 3.2, Theorem 4.4]{OzawaSuzuki} and \cite[Theorem B]{ESPJA} for the complete list of such characterizations.

\autoref{thm:TrAmSI} asserts that, for actions on sufficiently regular C*-algebras, the only difference between amenability and its tracial counterpart is given by equivariant property (SI), whose definition we recall below.

\begin{df}[{\cite[Proposition 5.1]{sato}, \cite[Definition 2.7]{GaborSI}}]\label{df-SI}
Let $G$ be a discrete group, let $A$ be a unital separable simple C*-algebra with $T(A)\neq \emptyset$, and let $\alpha\colon G\to\Aut(A)$ be an action. We say that $\alpha$ has \emph{equivariant property (SI)}\index{property (SI)} if for all positive contractions $e,f\in (A_\omega\cap A')^{\alpha_\omega}$ satisfying the following:
\begin{enumerate}
 \item $e$ is \emph{tracially small}, that is, $e\in J_{T(A)}$, and
\item $f$ is \emph{tracially supported at 1}, that is, 
\[\sup_{m\in \N} \|1-f^m\|_{2,T_\omega(A)}<1,\]
\end{enumerate}
there exists a contraction $s\in (A_\omega\cap A')^{\alpha_\omega}$ such that $s^*s=e$ and $fs=s$.
\end{df}

\begin{proof}[Proof of \autoref{thm:TrAmSI}]
Since all three conditions in the statement imply exactness of $G$, that is the existence of an amenable $G$-action on a compact metrizable space (see \cite[Corollary 6.2]{BusEchWil20} and \cite[Theorem 2.5]{ESPJA}), we may assume that $G$ is exact.

\eqref{SI1} $\Rightarrow$ \eqref{SI2}. The fact that all amenable actions are tracially amenable immediately follows from the definitions (see for instance the observation after \cite[Remark 2.3]{ESPJA}), and the main result of \cite{GaborSIamenable} shows that all amenable actions automatically have equivariant property (SI) (see also \cite[Lemma 5.3.3]{JulianPhD} for a proof of a special case).

\eqref{SI2} $\Rightarrow$ \eqref{SI3}. We will show that $(J_{T(A)} \cap A')^{\alpha_\omega}$ satisfies the Hjelmborg--R{\o}rdam criterion. Once we have done this, the conclusion will follow from \cite[Proposition 6.10]{CGSTW}. Let $e \in (J_{T(A)} \cap A')^{\alpha_\omega}_+$ be a nonzero positive contraction. Since $J_{T(A)}$ is an  \emph{equivariant $\sigma$-ideal} (see, for example, \cite[Proposition 7.3]{GHV}), there exists a positive contraction $x \in (J_{T(A)} \cap A')^{\alpha_\omega}$ such that $ex = e$. Set $f = 1 - x$ and note that $f$ is orthogonal to $e$. Moreover, $f$ is tracially supported at 1, hence by equivariant property (SI) there is a contraction $s\in (A_\omega\cap A')^{\alpha_\omega}$ such that $s^*s=e$ and $fs=s$, which implies $e ss^* = ef ss^* = 0$.

\eqref{SI3} $\Rightarrow$ \eqref{SI1}.
By \cite[Theorem 4.4]{OzawaSuzuki} it suffices to construct, for any action $\gamma\colon G\to \Aut(C)$ on a unital separable C*-algebra $C$, an equivariant unital completely positive map $\Phi\colon (C,\gamma) \to (A_\omega\cap A', \alpha_\omega)$.

By \cite[Theorem 2.5]{ESPJA}, tracial amenability implies the existence of an equivariant unital completely positive map
\[\varphi\colon (C,\gamma)\to (A^\omega\cap A',\alpha^\omega).\]
By a standard argument that goes back to \cite{KirchbergAbel}, and which we briefly spell out below for the reader's convenience, there exists an equivariant completely positive contractive map $\psi\colon (C,\gamma)\to (A_\omega\cap A',\alpha_\omega)$ making the following diagram commute: 
\[
\xymatrix{
&&A_\omega\cap A'\ar[d]^-{q}\\
C\ar[urr]^-{\psi}\ar[rr]_-{\varphi}	&&A^\omega\cap A'.
}
\]
Let indeed $\psi_0 \colon C \to A_\omega$ be a set-theoretic lift of $\phi$ to $A_\omega$, and let $D \subseteq A_\omega$ be a $G$-invariant subalgebra containing $\psi_0(C)$ and $A$. Arguing as in the proof of \eqref{SI2} $\Rightarrow$ \eqref{SI3}, since $J_{T(A)}$ is an {equivariant $\sigma$-ideal}, there is a positive contraction $b \in (J_{T(A)} \cap A')^{\alpha_\omega}$ such that $ba = a$ for all $a \in J_{T(A)} \cap D$. It is now straightforward to check that the function defined as $\psi(c) = (1 - b)\psi_0(c)$, for $c \in C$, is as desired.

%Use
%\cite[Proposition 1.18]{GaborLise}\footnote{Although explicitly recorded in \cite{GaborLise}, this elementary proof technique goes back to \cite{KirchbergAbel} and does not use all of the assumptions appearing in \cite{GaborLise}. The main idea is to choose a set-theoretical lift $s\colon C\to A_\omega$ and a $G$-invariant contraction $e\in (J_{T(A)}\cap A')^{\alpha_\omega}$ that acts as a unit on $s(C)$ (using a quasi-central quasi-invariant approximate unit). It can then be checked that $\psi(c)=(1-e)s(c)$ satisfies the desired properties.} to find an equivariant completely positive contractive map $\psi\colon (C,\gamma)\to (A_\omega\cap A',\alpha_\omega)$ making the following diagram commute: 

Note that $e=1-\psi(1)$ is a positive contraction in $(J_{T(A)}\cap A')^{\alpha_\omega}$, and let $\varepsilon>0$. By the Hjelmborg--R{\o}rdam criterion, there is $s \in (J_{T(A)}\cap A')^{\alpha_\omega}$ such that  $\| e - s^*s \| < \varepsilon$ and $\| (s^*s)(ss^*) \| < \varepsilon$. A standard saturation argument then yields $s \in (J_{T(A)}\cap A')^{\alpha_\omega}$ such that $e = s^*s$ and $e ss^* = 0$, which by the C*-identity gives $\psi(1) s = s$.

%and $f=\psi(1)$ are positive contractions in $(A_\omega\cap A')^{\alpha_\omega}$. Moreover $e\in J_{T(A)}$, so it is tracially small, while
%as $q(f)=\varphi(1)=1$, we have
%\[\sup_{m\in \N}\|1-f^m\|_{2,T^\omega(A)}=\sup_{m\in \N}\|q(1-f^m)\|_{2,T^\omega(A)}=\sup_{m\in \N}\|1-{q(f)}^m\|_{2,T^\omega(A)}=0,\]
%where $T^\omega(A)\subseteq T(A^\omega)$ denotes the induced set of limit traces.
%Hence $f$ is tracially supported at $1$. 

%Since $\alpha$ has equivariant property (SI), there exists a contraction $s\in (A_\omega\cap A')^{\alpha_\omega}$ satisfying $s^*s=1-\psi(1)$ and $\psi(1)s=s$.
Define $\Phi\colon C\to A_\omega\cap A'$ as
\[
\Phi(c)= \psi(c)+s^*\psi(c)s,\]
for all $c\in C$.
The map $\Phi$ is completely positive, being the sum of two completely positive maps. Furthermore 
\[\Phi(1)=\psi(1)+s^*\psi(1)s=\psi(1)+s^*s=1,\]
so $\Phi$ is unital.
Finally, $\Phi$ is equivariant since $s$ is $G$-invariant and $\psi$ is equivariant.
The result thus follows from \cite[Theorem 4.4]{OzawaSuzuki}.
\end{proof}
 
The following two examples give instances, one on a simple AF-algebra and the other on the Toeplitz algebra, of actions which are tracially amenable but not amenable. Using \autoref{thm:TrAmSI}, in \autoref{eg:SI} we are able to exhibit an equivariantly $\mathcal{Z}$-stable action for which equivariant property (SI) fails. %This is the first known example of an equivariantly $\mathcal{Z}$-stable
%action on a classifiable C*-algebra without equivariant 
%property (SI).

\begin{proof}[\autoref{eg:SI}]
	Denote by $X$ the Cantor set and fix a topologically amenable action $\gamma \colon F_n \curvearrowright X$ (let for instance $X=\partial F_n$ be the Gromov boundary of $F_n$ and consider the boundary action). 
	Using Elliott's classification of AF-algebras in \cite{Elliott}, it is possible to construct a unital simple AF-algebra $B$ whose trace space has extreme boundary $\partial_e T(B)$ homeomorphic to $X$, and
	an action $\beta\colon F_n\to \Aut(B)$ such that the induced action on $\partial_e T(B) \cong X$ is identified with $\gamma$ (see \cite[Example 2.11]{ESPJA} for the details). 
	
	Fix a pure state $\rho$ on $B$ and denote by $g_1,\dotsc,g_n$ the generators of $F_n$. 
	By \cite[Theorem 1.1]{KOS}, the action of all approximately inner automorphisms of $B$ on the pure state space of $B$ is transitive. Thus,
	there are approximately inner automorphisms $\varphi_1,\dotsc,\varphi_n\in \Aut(B)$ satisfying $\rho = \rho\circ \beta_{g_i} \circ \varphi_i$ for all $i=1,\dotsc,n$. 
	Denote by $\mathcal{Q}$ the universal UHF-algebra. Using the universal property of $F_n$, we define an action $\alpha$ of $F_n$ on $A= B\otimes \mathcal Q$ by setting
		\[\alpha_{g_i}= \beta_{g_i}\circ \varphi_i \otimes \id_\mathcal Q\]
		for all $i=1,\ldots,n$.
		
	Then $A$ is a unital simple AF-algebra and $\alpha$ is conjugate to $\alpha\otimes \id_\mathcal Q$. 
	As $\varphi_1,\dotsc,\varphi_n$ are approximately inner, they all act trivially on $T(B)$. Since moreover $\mathcal Q$ is monotracial, it follows that the action induced by $\alpha$ on $\partial_e T(A)\cong \partial_e T(B)\cong X$ is conjugate to that of $\beta$, and therefore can also be 
	identified with $\gamma$. It follows that the action on $T(A)$ induced by $\alpha$ is topologically amenable, which implies that $\alpha$ is tracially amenable by \cite[Theorem 2.5]{ESPJA}. On the other hand, the action induced by $\alpha$ on the state space $S(A)$
	of $A$ has a fixed point: indeed, given any state $\psi$
	on $\mathcal Q$, the state $\rho\otimes \psi$ is $F_n$-invariant. This means that the action on $S(A)$ is not topologically amenable, a necessary condition for amenability of $\alpha$ by \cite[Proposition 3.5]{OzawaSuzuki}. 
	\autoref{thm:TrAmSI} now implies that $\alpha$ does not have equivariant property (SI). 
	
	To see the failure of the weak containment property, let $\mathrm u \colon F_n \to \mathcal{B}(H)$ be the universal representation of $F_n$. Identify $\rho \otimes \psi$
	with the map sending $a \in A$ to $(\rho\otimes \psi)(a) \cdot 1_{\mathcal{B}(H)}$. Since $\rho \otimes \psi$ is $F_n$-invariant, we get that the pair $(\rho \otimes \psi, \mathrm u)$
	satisfies, for all $a \in A$ and $g \in F_n$, the covariance relation
	\[
	(\rho \otimes \psi)(\alpha_g(a))\cdot 1_{\mathcal{B}(H)} = (\rho \otimes \psi)(a)  \cdot 1_{\mathcal{B}(H)} = \mathrm u_g (\rho \otimes \psi)(a)  \cdot 1_{\mathcal{B}(H)} \mathrm u_{g^{-1}}.
	\]
	By \cite[Lemma 5.3]{BusEchWil22} we therefore get that the function
	\begin{align*}
	(\rho \otimes \psi) \rtimes \mathrm u \colon C_c(F_n, A) &\to C^\ast_{\mathrm{max}}(F_n) \\
	f &\mapsto \sum_{g \in F_n} (\rho \otimes \psi)(f(g)) \mathrm{u}_g
	\end{align*}
	extends to a completely positive map $(\rho \otimes \psi) \rtimes \mathrm u \colon A\rtimes_{\mathrm{max}} F_n \to C^\ast_{\mathrm{max}}(F_n)$. If $\alpha$ had the weak containment property,
	then $A\rtimes_{\mathrm{max}} F_n \cong A\rtimes_{\mathrm{red}} F_n$, with the isomorphism being the canonical quotient map. The restriction of such map to $C^\ast_{\mathrm{red}}(F_n)
	\subseteq A\rtimes_{\mathrm{red}} F_n$ would then be an isomorphism between $C^\ast_{\mathrm{red}}(F_n)$ and $C^\ast_{\mathrm{max}}(F_n)$ extending the identity map on $C_c(F_n)$, a contradiction
	since $F_n$ is not amenable.
	\end{proof}

\begin{proof}[\autoref{eg-Toeplitz}]
	Denote by $\mathbb{T}\subseteq \C$ the unit circle and by $H^2\subseteq L^2(\mathbb{T})$ the Hardy space, that is, the Hilbert space generated by all the monomials $z^n\in L^2(\mathbb{T})$ for $n \ge 0$. Let $P\colon L^2(\mathbb{T})\to H^2$ be the orthogonal projection onto the Hardy space. For a function $f\in C(\mathbb{T})$, let $M_f\colon L^2(\mathbb{T})\to L^2(\mathbb{T})$ be the corresponding multiplication operator, and set $T_f=PM_fP\in \mathcal B(H^2)$. The C*-algebra $\mathcal{T}$ generated by $\{T_f\colon f\in C(\mathbb{T})\}$ is the Toeplitz algebra, and there is a short exact sequence 
	\[
\xymatrix{		0\ar[r]& \mathcal K(H^2)\ar[r]& \mathcal{T}\ar[r]^-{\pi}& C(\mathbb{T})\ar[r]& 0},
	\]
	where the quotient map $\pi$ is determined by $\pi(T_f)= f$, for all $f\in C(\mathbb{T})$; see for instance \cite[Theorem 3.5.11]{murphy}.
	
Let $\{e_1,e_2\}$ denote the canonical orthonormal basis 
of $\mathbb{R}^2$, and
denote by $\mathbb{R} \mathrm P^1$ the space of one-dimensional subspaces of $\mathbb{R}^2$. It is easy to see that the natural action of the special linear group $\mathrm{SL}(2,\mathbb{R})$ on $\mathbb{R} \mathrm P^1$ is transitive, and that the stabilizer of $\mathbb{R}\cdot e_1$ is given by the subgroup $P\subseteq \mathrm{SL}(2,\mathbb{R})$ of upper triangular matrices. Note that $P$ is solvable and therefore amenable. Since $F_n$ embeds into $\mathrm{SL}(2,\mathbb{R})$ as a discrete subgroup\footnote{This is standard: for example, one can use the ping-pong lemma to see that the subgroup of $\mathrm{SL}(2,\mathbb{R})$ generated by $\begin{psmallmatrix}1 & 2\\0 & 1\end{psmallmatrix}$ and $\begin{psmallmatrix}1 & 0\\2 & 1\end{psmallmatrix}$ is isomorphic to $F_2$.}, the restriction $\gamma$ of $\mathrm{SL}(2,\mathbb{R}) \curvearrowright \mathbb{R} P^1$ to $F_n$ is thus amenable by \cite[Examples 2.2.18(1)]{ADRenault2000amenableGroupoids}.
Using that $\mathbb{R} \mathrm P^1$ is homeomorphic to $\mathbb{T}$, we obtain a topologically amenable action of $\mathrm{SL}(2,\mathbb{R})$ on $\mathbb{T}$ that restricts to a topologically amenable action of $F_n$ on $\mathbb{T}$. 
	
	To each $g\in F_n$ we now associate a Toeplitz operator $T_g$ by identifying $g$ with an element of $C(\mathbb{T})$ via its action on $\mathbb{T}\subseteq \C$. 
	Since $\mathrm{SL}(2,\mathbb{R})$ is path-connected, every $g\in F_n$ acts on $\mathbb{T}$ via a homeomorphism that is homotopic to the identity and thus with winding number $1$. By the Toeplitz Index Theorem (see, for example, \cite[Theorem 2.3.2]{AnalyticKhomology}) this implies that the Toeplitz operators $T_g\in \mathcal{T}$ have Fredholm index equal to $-1$. By Brown--Douglas--Fillmore theory (\cite[Theorem 3.1]{BDF}) we conclude that there exist compact operators $K_1,\dotsc,K_n\in \mathcal K(H^2)$ and unitaries $U_1,\dotsc,U_n\in \mathcal U(H^2)$ such that 
	\[T_{g_i} + K_i = U_i T_z U_i^*\]
	for $i=1,\ldots,n$, 
	where $z\in C(\mathbb{T})$ is the identity function.
	
	Define finally an action $\alpha$ of $F_n$ on $\mathcal{T}$ by letting each generator $g_i$ act by $\Ad(U_i)$. 
	By construction, $\alpha$ induces the given topologically amenable action $\gamma$ on the trace space $T(\mathcal{T}) \cong\Prob(\mathbb{T})$ and is therefore tracially amenable by \cite[Theorem 2.5]{ESPJA}. However, if $\alpha$ were amenable, so would be its restriction to the invariant ideal $\mathcal K(H^2)$ by \cite[Proposition 3.23]{BusEchWil22}. Since amenability is preserved by Morita equivalence (\cite[Proposition 3.20]{BusEchWil22}), this would imply that the trivial action $F_n\curvearrowright \mathbb{C}$ is amenable, which is a contradiction. We conclude that $\alpha$ is not amenable, as desired.
	\end{proof}

We conclude the paper with the proof of \autoref{thm-noembedding}. Given a discrete group $G$, a \emph{$G$-boundary} is a compact Hausdorff space $X$ along with an action $G \curvearrowright X$ which is \emph{minimal}, meaning that all orbits are dense, and \emph{strongly proximal}, which is the case if the weak$^*$-closure of the orbit of every probability measure $\mu \in \text{Prob}(X)$ contains a Dirac measure.

\begin{proof}[Proof of \autoref{thm-noembedding}]
	Assume by contradiction that $\Phi$ exists. Let $Y\subseteq T(A)$ be a minimal closed $G$-invariant subspace and denote by $\mathrm{ev}\colon A\to C(Y)$ the evaluation map, that is $\mathrm{ev}(a)(\tau)=
	\tau(a)$ for $a \in A$ and $\tau \in Y$.
	The composition 
	\[C(X)\xrightarrow{\Phi}A\xrightarrow{\mathrm{ev}}C(Y)\]
	is a unital positive equivariant map, hence by  \cite[Lemma 3.10]{Kal_Ken_2017} it
	is a unital injective $*$-homomorphism  and it  is therefore induced by a continuous surjective map $\phi\colon Y\to X$.
	
	Fix $x\in X$, a non-zero positive function $f\in C(X)$ with $f(x)=0$, and $\tau\in Y$ with $\phi(\tau)=x$.  
	Then the element $(\ev\circ \Phi)(f)\in C(Y)$ is non-zero and satisfies 
		\[\tau(\Phi(f))=(\ev\circ\Phi)(f)(\tau)=f(\phi(\tau))=0.\]
Since $A$ is simple, the trace $\tau$ is faithful. As $\Phi(f)$ is positive,
we deduce from the above that $\Phi(f)=0$, contradicting the fact that $(\mathrm{ev}\circ \Phi)(f)\not=0$. This contradiction shows that $\Phi$ does not exist and proves the proposition.
\end{proof}

% \bibliography{references}{}
%	\bibliographystyle{alpha}

\newcommand{\etalchar}[1]{$^{#1}$}

\end{document}